\documentclass[fleqn]{mat01}
\usepackage{times,mathtimy,amssymb,latexsym}
\begin{document}

\setcounter{page}{309}
\firstpage{309}

\def\d{\mbox{\rm d}}

\newtheorem{theore}{Theorem}
\renewcommand\thetheore{\arabic{section}.\arabic{theore}}
\newtheorem{theor}[theore]{\bf Theorem}
\newtheorem{lem}[theore]{Lemma}
\newtheorem{rem}[theore]{Remark}
\newtheorem{coro}[theore]{\rm COROLLARY}

\newtheorem{step}{Step}

\renewcommand\theequation{\arabic{section}.\arabic{equation}}

\title{Estimates and nonexistence of solutions of the scalar curvature
equation on noncompact manifolds}

\markboth{Zhang Zonglao}{Scalar curvature equation on noncompact manifolds}

\author{ZHANG ZONGLAO}

\address{Department of Mathematics, Wenzhou Normal College, Wenzhou,
Zhejiang 325 000, People's Republic of China\\
\noindent E-mail: zonglao@sohu.com}

\volume{115}

\mon{August}

\parts{3}

\pubyear{2005}

\Date{MS received 7 January 2004}

\begin{abstract}
This paper is to study the conformal scalar curvature equation on
complete noncompact Riemannian manifold of nonpositive curvature. We
derive some estimates and properties of supersolutions of the scalar
curvature equation, and obtain some nonexistence results for complete
solutions of scalar curvature equation.
\end{abstract}


\maketitle

\section{Introduction}

Let $M$ be an $n$-dimensional complete Riemannian manifold with metric
$g_{0}$ and scalar curvature $k (x)$. The problem of the conformal
deformation of scalar curvature is to find conditions on the function $K
(x)$ so that $K (x)$ is the scalar curvature of a conformally related
metric $g = \rho (x) g_{0}$, here $\rho (x)$ is some positive function
on $M$. As is well-known, let $\rho = u^{4/(n - 2)}$ for $n \geq 3$.
Then this problem is equivalent to finding a positive solution of the
following equation (the so-called (conformal) scalar curvature
equation):
\begin{equation}
c_{n} \Delta u - ku + Ku^{\sigma} = 0,
\end{equation}
where $c_{n} = {\frac{4 (n - 1)}{n - 2}}, \sigma = {\frac{n + 2}{n -
2}}, \Delta u$ is the Laplacian of $u$ with respect to the metric
$g_{0}$. If $u$ is a positive solution of eq.~(1.1) such that the metric
$g = u^{4/(n - 2)} g_{0}$ is complete, we call $u$ to be complete. If $M$
is noncompact, usually one seeks complete solutions of eq.~(1.1).

The problem of conformal deformation of scalar curvature has been
extensively studied by many authors in recent years \cite{1,2,7,8,9}.
However, this problem is far from being settled, especially if $M$ is
noncompact.

The purpose of this paper is to continue to study the problem of
conformal deformation of scalar curvature. We consider the case when $M$
is complete and noncompact. We will derive some estimates and properties
of supersolutions of eq.~(1.1), and obtain some nonexistence results for
complete solutions of (1.1). If $M$ is noncompact, in order to obtain a
positive solution of (1.1), one usually exploits the method of
supersolution--subsolution. To do so, one usually needs a positive
supersolution bounded from below by a positive constant or a positive
subsolution bounded from above by a constant as in \cite{1,2,7,8,9}. Here
we will see that, under some suitable assumptions, the supersolutions
of (1.1) have no positive constant lower bound for a large class of the
functions $K$.

This paper is organized as follows. In \S2, we introduce some notations
and give some results that will be needed in \S3. In \S3 we will state
and prove the main results of our paper.

\section{Preliminaries}

We call $M$ a CH manifold if it is a complete simply-connected $C^{\infty}$
Riemannian manifold of nonpositive sectional curvature \cite{6}. If $M$ is a
CH manifold, by the well-known Cartan--Hadamard Theorem, for every $o \in
M$, the exponential mapping $\exp_{o}$: $M_{o} \rightarrow M$ is a
diffeomorphism \cite{3}, where $M_{o}$ denotes the tangent space to $M$
at $o$. This diffeomorphism gives a global normal coordinate neighborhood of $M$
center at $o$. Denote by $x$ the coordinates of points and by $(r,
\theta)$ the (geodesic) polar coordinates around $o$, where $r = r (x)
\equiv \ \hbox{dist} \ (o, x)$ is the geodesic distance from $o$.

We call a CH manifold $M$ to be strongly symmetric around $o \in M$ if
every linear isometry $\phi$: $M_{o} \rightarrow M_{o}$ is realized as
the differential of an isometry $\Phi$: $M \rightarrow M$, i.e., $\Phi
(o) = o$ and $\Phi_{*} (o) = \phi$, where $\Phi_{*} (o)$ denotes the differential of $\Phi$ at $o$.
For a more detailed discussion about strongly symmetric manifold, we
refer the reader to \cite{6} (where the authors use the term `model'
instead of `strongly symmetric manifold').

From now on we will assume that $M$ is an $n$-dimensional strongly
symmetric CH manifold around $o$, where $o$ is a fixed point in $M$. Let
$g_{0}$ be the metric of $M$ and $k (x)$ the scalar curvature of
$g_{0}$. We always assume $n = \dim M \geq 3$.

In the polar coordinates, the metric $g_{0}$ is expressed by
\setcounter{equation}{0}
\begin{equation}
\d s^{2} = \d r^{2} + \sum\limits_{i, j} \d_{ij} \d \theta^{i} \d
\theta^{j} = \d r^{2} + h (r)^{2} \d \Theta^{2}
\end{equation}
on $M - \{ o \}$ \cite{6}. Here $\d_{ij} = g_{0} \big(
{\frac{\partial}{\partial \theta^{i}}}, {\frac{\partial}{\partial
\theta^{j}}} \big)$ and $\d \Theta^{2}$ denotes the canonical metric
on the unit sphere of $M_{o}$. Let $S_{r}$ be the geodesic sphere of
$M$ with center $o$ and radius $r$. The Riemannian volume element of
$S_{r}$ can be written as
\begin{equation}
\d S_{r} = \sqrt{D (r, \theta)} \d \theta^{1} \cdots \d \theta^{n - 1},
\end{equation}
where $D \equiv \det (\d_{ij})$. We will denote by $V (r)$ the volume of
$S_{r}$.

If $u (r)$ is a $C^{2}$ function defined on $(0, \infty)$, we can
consider it as a function defined on $M - \{ o \}$. A calculation shows
\begin{align}
\Delta u &= {\frac{1}{\sqrt{D}}} \partial_r (\sqrt{D} \partial_{r} u),\\[.3pc]
\Delta r &= {\frac{1}{\sqrt{D}}} \partial_r (\sqrt{D}) = \partial_{r} \log \sqrt{D}
\end{align}
and
\begin{equation}
\Delta u = u'' + (\Delta r) u'.
\end{equation}

We can define a scalar product operation $\eta$ on $M$ as follows:
\begin{equation}
\eta : \textbf{R} \times M \rightarrow M; \ (t, (r, \theta)) \mapsto
(tr, \theta).
\end{equation}

We also write $\eta (r, x) \equiv rx$.

\section{The main results}

In this section we will state and prove our main results. To do so, we
first introduce a notation.

Let $M$ be an $n$-dimensional strongly symmetric CH manifold around $o$.
If $f$ is a continuous function on $M$, define
\setcounter{equation}{0}
\begin{equation}
\bar{f} (r) \equiv {\frac{1}{V (r)}} \int_{S_{r}} f \d S_{r} =
{\frac{1}{V (1)}} \int_{S_{1}} [f (r \xi)] \d S_{1}, \quad (\xi \in
S_{1})
\end{equation}
(the second equality is by (2.1)).

\begin{theor}[\!]
Let $M$ be an $n$-dimensional strongly symmetric CH manifold around $o$
with metric $g_{0}$. Let $k$ be the scalar curvature of the metric
$g_{0}${\rm ,} and let $K \in C^{\infty} (M)$. Suppose $n \geq 3$. If $u$ is a
$C^{2}$ positive supersolution of equation~$(1.1)$ on $M${\rm ,} set $\alpha = 1 -
\sigma$ and $v = u^{\alpha}${\rm ,} then $\bar{v}$ satisfies the
following inequality
\begin{equation}
\bar{v} (r) \geq {\frac{1}{n - 1}} \int_{0}^{r} {\frac{1}{V
(s)}} \int_{B (s)} K \d \mu \d s
\end{equation}
for all $r \geq 0${\rm ,} here $B (s)$ is the geodesic ball of radius $s$ and
center $o${\rm ,} and $\d \mu$ the volume element of $M$.
\end{theor}

\begin{proof}$\left.\right.$

\begin{step}{\rm We first prove that $v$ satisfies the
following inequality
\begin{equation}
\Delta v \geq {\frac{1}{n - 1}} (K - kv).
\end{equation}

In fact, a computation shows that
\begin{equation}
\Delta v = \alpha u^{\alpha - 1} \Delta u + {\frac{\alpha - 1}{\alpha}}
u^{- \alpha} |\Delta u^{\alpha}|^{2}.
\end{equation}
Since $u$ is a positive supersolution of eq.~(1.1) and $\alpha = -
{\frac{4}{n - 2}} < 0$, we have
\begin{align*}
\Delta v &\geq {\frac{\alpha}{c_{n}}} u^{\alpha - 1} (ku -
Ku^{\sigma})\\[.3pc]
&= {\frac{1}{n - 1}} (K - kv).
\end{align*}}
\end{step}

\begin{step}{\rm
We prove that $\bar{v}$ satisfies
\begin{equation}
\Delta \bar{v} \geq {\frac{1}{n - 1}} (\bar{K} - k
\bar{v}).
\end{equation}
By the definition (3.1) of $\bar{v}(r)$, for any $r > 0$, we have
\begin{equation*}
\bar{v}'(r) = \frac{1}{V(1)} \int_{S_{1}} \frac{\partial
v}{\partial r} (r \xi) \d S_{1} = \frac{1}{V(r)} \int_{S_{r}}
\partial_{r} v \d S_{r}.
\end{equation*}
So for $r > 0$,
\begin{equation}
\bar{v}' (r) V(r) = \int_{S_{r}} \partial_{r} v \d S_{r}.
\end{equation}
Thus, for $r > 0$, by the divergence theorem \cite{4,5} we have
\begin{equation}
\int_{B(r)} \Delta v \d \mu = \int_{S_{r}} \partial_{r} v \d S_{r} =
\bar{v}'(r) V(r).
\end{equation}
But we also have
\begin{equation}
\int_{B(r)} \Delta v \d \mu = \int_{0}^{r} \int_{S_{t}} \Delta v \d
S_{t} \ \d t.
\end{equation}
So from (3.7) and (3.8), we obtain
\begin{align}
\int_{S_{r}} \Delta v \d S_{r} &= \left[\bar{v}' (r) V(r) \right]'\nonumber\\[.3pc]
&= V(r) \left\lbrace \bar{v}'' (r) + \frac{V'(r)}{V(r)}
\bar{v}' (r) \right\rbrace.
\end{align}

On the other hand, by (2.2) and (2.4),
\begin{equation*}
V'(r) = \int_{0}^{2 \pi} \ \int_{0}^{\pi} \cdots \int_{0}^{\pi}
\frac{\partial_{r} \sqrt{D}}{\sqrt{D}} \sqrt{D} \d \theta^{1} \cdots
\d \theta^{n - 2} \d \theta^{n - 1} = V(r) \Delta r,
\end{equation*}
so we get
\begin{equation}
\Delta r = \frac{V'(r)}{V(r)}.
\end{equation}

From (2.5), (3.9) and (3.10) we get
\begin{align*}
\int_{S_{r}} \Delta v \d S_{r} &= V(r) \left\lbrace \bar{v}'' (r)
+ (\Delta r) \bar{v}' (r) \right\rbrace\\[.3pc]
&= V(r) \Delta \bar{v} (r).
\end{align*}
Thus, for $r > 0$, we obtain
\begin{equation}
\Delta \bar{v} (r) = \frac{1}{V(r)} \int_{S_{r}} \Delta v \d S_{r}.
\end{equation}

Now from (3.3) and (3.11) we see that
\begin{align*}
\Delta \bar{v} (r) &\geq \frac{1}{(n - 1) V(r)} \int_{S_{r}} (K -
kv) \d S_{r}\\[.3pc]
&= \frac{1}{n - 1} (\bar{K} - k \bar{v}).
\end{align*}}
\end{step}

\begin{step}{\rm
We are now ready to prove inequality (3.2).

Integrating (3.5) we get
\begin{align}
\int_{B(r)} \Delta \bar{v} (r) \d \mu &\geq \frac{1}{n - 1}
\int_{B(r)} (\bar{K} - k \bar{v})\d \mu\nonumber\\[.3pc]
&= \frac{1}{n - 1} \int_{0}^{r} \{ \bar{K}
(t) - k(t) \bar{v} (t) \} V(t) \d t.
\end{align}

It is easy to prove that $\bar{v}' (r)$ is continuous on $[0,
\infty)$ and $\bar{v}' (0) = 0$. Then, by the divergence theorem,
(3.12) implies that
\begin{align*}
\bar{v}' (r) V(r) &= \int_{S_{r}} \bar{v}' (r) \d S_{r}\\[.3pc]
&= \int_{B(r)} \Delta \bar{v} (r) \d \mu\\[.3pc]
&\geq \frac{1}{n - 1} \int_{0}^{r} \{ \bar{K} (t)
- k(t) \bar{v} (t) \} V(t) \d t,
\end{align*}
that is,
\begin{equation}
\bar{v}' (r) \geq \frac{1}{(n  -1) V(r)} \int_{0}^{r} \{
\bar{K} (t) - k(t) \bar{v} (t) \} V(t) \d t.
\end{equation}
Note that $k \leq 0$. Integrating (3.13) we obtain
\begin{align}
\bar{v} (r) &\geq \bar{v} (0) + \frac{1}{n - 1} \int_{0}^{r}
\frac{1}{V(s)} \int_{0}^{s} \{ \bar{K} (t) - k(t)
\bar{v} (t) \} V(t)\, \d t \d s\nonumber\\[.4pc]
&\geq \frac{1}{n - 1} \int_{0}^{r} \frac{1}{V (s)} \int_{0}^{s}
\bar{K} (t) V(t)\, \d t \d s\nonumber\\[.4pc]
&= \frac{1}{n - 1} \int_{0}^{r} \frac{1}{V (s)} \int_{B(s)} K \d \mu \d
s.
\end{align}}
\end{step}
Thus (3.2) is established, and this completes the proof of Theorem~3.1.
\hfill $\Box$
\end{proof}

\begin{theor}[\!]
Let $M$ be an $n$-dimensional strongly symmetric CH manifold around $o$
with $n \geq 3$. Let $k$ be the scalar curvature of $M${\rm ,} and let $K \in
C^{\infty} (M)$. Assume $u$ is a $C^{2}$ positive supersolution of
eq.~$(1.1)$ on $M$. If we have either
\begin{enumerate}
\renewcommand\labelenumi{{\rm (\alph{enumi})}}
\item
\begin{equation}
\int_{0}^{\infty} \frac{1}{V (s)} \int_{B(s)} K {\rm d} \mu {\rm d} s = + \infty
\end{equation}
or

\item for $r$ large{\rm ,} $\int_{B(r)} K {\rm d} \mu \geq 0${\rm ,} and
\begin{equation}
\int_{0}^{\infty} \frac{1}{V(s)} \int_{B(s)} |k| {\rm d} \mu {\rm d} s = +
\infty{\rm ,}
\end{equation}
then
\begin{equation}
\inf\limits_{x \in M} u (x) = 0.
\end{equation}
\end{enumerate}
\end{theor}

\begin{proof}
As in Theorem~3.1, set $\alpha = 1 - \sigma$ and $v = u^{\alpha}$, then
$\bar{v}$ satisfies inequalities (3.13) and (3.14).

If (a) holds, it is obvious that $\bar{v} (\infty) = \infty$, and
hence $\sup_{x \in M} v(x) = \infty$. Since $\alpha < 0$, we see that
$\inf_{x \in M} u (x) = 0$.

If (b) holds, then there exists $\tau > 0$ such that for all $r \geq
\tau, \int_{B(r)} K \d \mu \geq 0$. From (3.13) we have that for $r \geq
\tau, \bar{v}' (r) \geq 0$ and hence $\bar{v} (r)$ is
increasing. This means that if we set $C = \inf \bar{v} (r)$, then
$C > 0$. From (3.14) we have
\begin{align*}
\bar{v} (r) &\geq \frac{1}{n - 1} \int_{0}^{r} \frac{1}{V (s)}
\int_{0}^{s} \bar{K} (t) V(t) \d t \d s\\[.5pc]
&\quad\, + \frac{C}{n - 1}
\int_{0}^{r} \frac{1}{V (s)} \int_{0}^{s} |k(t)| V (t) \d t \d s\\[.5pc]
&= \frac{1}{n - 1} \int_{0}^{r} \frac{1}{V (s)} \int_{B(s)} K \d \mu \d
s + \frac{C}{n - 1} \int_{0}^{r} \frac{1}{V (s)} \int_{B(s)} |k| \d \mu
\d s\\[.5pc]
&= \frac{1}{n - 1} \int_{0}^{\tau} \frac{1}{V (s)} \int_{B(s)} K \d \mu
\d s + \frac{1}{n - 1} \int_{\tau}^{r} \frac{1}{V (s)} \int_{B(s)} K \d \mu
\d s\\[.5pc]
&\quad + \frac{C}{n - 1} \int_{0}^{r} \frac{1}{V (s)} \int_{B(s)} |k|
\d \mu \d s.
\end{align*}
Now it is easy to see that $\bar{v} (\infty) = \infty$. Then (3.17)
follows as in Case (a). The proof of Theorem~3.2 is finished. \hfill
$\Box$
\end{proof}

To show the following theorems, we need a lemma. This result itself is
an interesting property of CH manifolds.

\begin{lem}

Let $M$ be an $n$-dimensional strongly symmetric CH manifold around $o$.
As before{\rm ,} let $S_{r}$ be the geodesic sphere of $M$ with center $o$ and
radius $r${\rm ,} and $V(r)$ the volume of $S_{r}$. Then for every $\delta \in
[0, 1)${\rm ,}
\begin{equation}
\lim_{r \rightarrow \infty} \frac{V(r)}{r^{n - 2 + \delta}} = + \infty.
\end{equation}
\end{lem}

\begin{proof}
Since the sectional curvatures of $M$ are less than or equal to zero, by
the volume comparison theorem (comparing with the $n$-dimensional
Euclidean space ${\bf R}^{n}$) \cite{3,4}, we have vol$(B(r)) \geq
r^{n} \omega_{n}$, where vol$(B(r))$ denotes the volume of the ball
$B(r) = \{ x \in M, \ \hbox{dist}(o, x) < r \}$ and $\omega_{n}$
denotes the volume of the unit ball in ${\bf R}^{n}$. Hence we\break get
\begin{equation}
\lim\limits_{r \rightarrow \infty} \frac{\int_{0}^{r} V(t) \d t}{r^{n -
1 + \delta}} = \lim\limits_{r \rightarrow \infty} \frac{\hbox{vol} (B(r))}{r^{n -
1 + \delta}} = + \infty.
\end{equation}

On the other hand, by the Laplacian comparison theorem (see p.~26 of
\cite{6}) (again comparing with the $n$-dimensional Euclidean space
${\bf R}^{n}$), we have
\begin{equation}
\Delta r \geq \Delta_{0} r = \frac{n - 1}{r},
\end{equation}
for $r > 0$. Here $\Delta_{0}$ denotes the Laplacian of ${\bf R}^{n}$. So for
$r > 0$ we have
\begin{align*}
\left(\frac{V(r)}{r^{n - 2 + \delta}}\right)' &= \frac{V(r)}{r^{n - 2 +
\delta}} \left(\frac{V'(r)}{V(r)} - \frac{n - 2 + \delta}{r}\right)\\[.3pc]
&= \frac{V(r)}{r^{n - 2 +
\delta}} \left(\Delta r - \frac{n - 2 + \delta}{r}\right)\\[.3pc]
&> 0.
\end{align*}
This means that $\frac{V(r)}{r^{n - 2 + \delta}}$ is increasing and
hence $\lim_{r \rightarrow \infty} \frac{V(r)}{r^{n - 2 + \delta}}$
exists (may be $+ \infty$). Then using the L' H$\hat{\hbox{o}}$pital's rule, we
have
\begin{equation}
\lim\limits_{r \rightarrow \infty} \frac{\int_{0}^{r} V(t) \d t}{r^{n -
1 + \delta}} = \lim\limits_{r \rightarrow \infty} \frac{V(r)}{(n - 1 +
\delta) r^{n - 2 + \delta}}.
\end{equation}
Then (3.18) follows immediately from (3.19) and (3.21).
\hfill $\Box$
\end{proof}

\begin{theor}[\!]
Let $M$ be an $n$-dimensional strongly symmetric CH manifold around $o$
with $n \geq 3$. Let $k$ be the scalar curvature of $M${\rm ,} and let $K \in
C^{\infty} (M)$ such that{\rm ,} for $r$ large{\rm ,} $\int_{B(r)} K {\rm d} \mu \geq 0$.
Assume $u$ is a $C^{2}$ positive supersolution of eq.~$(1.1)$ on $M$. If
we have either
\begin{enumerate}
\renewcommand\labelenumi{\rm (\alph{enumi})}
\item $\lim_{r \rightarrow \infty} \int_{B(r)} K {\rm d} \mu = + \infty${\rm ,} and
\begin{equation}
\hskip -1.25pc \lim\limits_{r \rightarrow \infty} \frac{r^{2} \bar{K} (r)}{r
\Delta r - 1} = \beta, \quad \hbox{where} \quad \beta > 0 \ \
\hbox{or}\ \ \beta = + \infty
\end{equation}
or
\item $\lim_{r \rightarrow \infty} \int_{B(r)} |k| {\rm d} \mu = + \infty$, and
\begin{equation}
\hskip -1.25pc \lim\limits_{r \rightarrow \infty} \frac{r^{2} |k (r)|}{r
\Delta r - 1} = \gamma, \quad \hbox{where} \quad \gamma > 0 \
\ \hbox{or}\ \ \gamma = + \infty{\rm ,}
\end{equation}
then
\begin{equation}
\hskip -1.25pc \inf_{x\in M} u (x) = 0.
\end{equation}
\end{enumerate}
\end{theor}

\begin{proof}
We only prove that condition (a) implies (3.24). The proof for (b)
is similar to that of (a). From Theorem 3.2, by the limit
comparison test for improper integrals, we only need to prove that
$\lim_{r\rightarrow \infty}\frac{r}{V(r)}\int_{B(r)}K\d\mu =
\lambda$, where $\lambda > 0$ or $\lambda = + \infty$. In fact,
from the assumption of the theorem and Lemma (3.3) we have
\begin{equation*}
\lim\limits_{r\rightarrow\infty} \int_{0}^{r}\bar{K}(t) V(t) \d t
= \lim\limits_{r\rightarrow\infty} \frac{V(r)}{r} = + \infty.
\end{equation*}

So by the L'H$\hat{\rm o}$pital's rule, we have
\begin{align*}
\lim\limits_{r\rightarrow\infty} \frac{r}{V(r)} \int_{B(r)} K \d
\mu  &= \lim\limits_{r\rightarrow\infty}
\frac{\int_{0}^{r}\bar{K}(t)V(t)\d t}{\frac{V(r)}{r}}\\[.3pc]
&= \lim\limits_{r\rightarrow\infty}\frac{\bar{K}(r)V(r)}
{\left(\frac{V(r)}{r}\right)'}
\end{align*}
\begin{align*}
&= \lim\limits_{r\rightarrow\infty}\frac{r^{2}\bar{K}(r)}{r\Delta
r - 1}\\[.3pc]
&=\beta.
\end{align*}
This completes the proof of the theorem. \hfill $\Box$
\end{proof}

\begin{rem}{\rm
Let ${\bf H}^{n} (-c^{2})$ be the $n$-dimensional
(simply-connected and complete) hyperbolic space form with
constant sectional curvature $-c^{2}$ $(c > 0)$. If $M = {\bf H}^{n}
(-c^{2})$ with $n \geq 3$, then for any $K \in
C^{\infty}(M)$ such that $\int_{B(r)} K\d \mu \geq 0$ for $r$
large, it is easy to verify that the condition (b) of Theorem~3.4
is satisfied since $\Delta r = (n - 1)c\coth (cr)$, and hence for
any $C^{2}$ positive supersolution $u$ of eq.~(1.1), we have
$\inf_{x\in M}u(x) = 0$. Comparing  with theorem~1.4 in
\cite{8}, this is very different from the case of the Euclidean
space. Hence the method used by Ni to obtain the existence results
for eq.~(1.1) in ${\bf R}^{n}$ is no longer valid in the
case of hyperbolic space forms.}
\end{rem}

\begin{theor}[\!]
Let $M$ be a $CH$ manifold with metric $g_{0}$. Suppose $M$ is
strongly symmetric around $o$ with respect to the metric $g_{0}$
and $n = \dim (M) \geq 3$. Let $k(r)$ be the scalar curvature of
the metric $g_{0}$ and $K(r(x)) \in C^{\infty}(M)$. If there is $a
> 0$ such that{\rm ,} for $r \geq a,
\int_{0}^{r}\frac{1}{V(t)}\int_{B(t)} K\d \mu\d t > 0$ and
\begin{equation}
\int_{a}^{\infty} \left[\int_{0}^{r}\frac{1}{V(t)}\int_{B(t)} K\d
\mu \d t \right]^{-1/2} \d r < + \infty,
\end{equation}
then there exists no complete metric $g$ on $M$ such that

\begin{enumerate}
\renewcommand\labelenumi{\rm (\alph{enumi})}
\item $g$ is {\rm (}pointwise{\rm )} conformal to the metric $g_{0}${\rm ,}
\item $M$ is strongly symmetric around $o$ with respect to the
metric $g${\rm ,}
\item $K(r)$ is the scalar curvature of the metric $g$.
\end{enumerate}
\end{theor}

\begin{proof}
We argue by contradiction. Assume that $g$ is a complete metric on
$M$ satisfying (a)--(c). Since $g$ is conformal to the metric
$g_{0}$, there exists a $C^{\infty}$ positive function $u$ such
that $g = u^{4/(n -2)}g_{0}$ and $u$ is a solution of eq.~(1.1).
It is obvious that $u = u(r)$ since the metrics $g$ and $g_{0}$
are both strongly symmetric around $o$. Set $v(r) = [u(r)]^{-4/(n
-2)}$. Then by Theorem~(3.1), $v(r)$ satisfies
\begin{equation*}
v(r) \geq \frac{1}{n - 1}\int_{0}^{r} \frac{1}{V(t)} \int_{B(t)}
K\d \mu\d t.
\end{equation*}

Therefore we get
\begin{equation*}
[u(r)]^{2/(n - 2)}\leq \left\lbrace\frac{1}{n -
1}\int_{0}^{r}\frac{1}{V(t)}\int_{B(t)}K\d \mu\d t
\right\rbrace^{-1/2}
\end{equation*}
for $r \geq a$. By (3.25), this implies
\begin{align*}
\int_{a}^{\infty}[u(r)]^{2/(n - 2)}\d r \leq\int_{a}^{\infty}
\left\lbrace\frac{1}{n -
1}\int_{0}^{r}\frac{1}{V(t)}\int_{B(t)}K\d \mu \d
t\right\rbrace^{-1/2}\d r < + \infty.
\end{align*}

That is, for the metric $g$, the ray $\varphi = \{(r,
\theta_{0})|a \leq r < \infty\}$ has finite length for a fixed
$\theta_{0}$. Thus the metric $g$ is not complete. This is a
contradiction. \hfill $\Box$
\end{proof}

\begin{coro}$\left.\right.$\vspace{.5pc}

\noindent Let $M, g_{0}, o$ and $k(r)$ be as in Theorem~$3.6$. Assume the
sectional curvature $\sec (g_{0})$ of $g_{0}$ satisfies
$-c^{2} \leq \sec(g_{0}) \leq 0$ for some positive constant
$c$. If $K(r) \in C^{\infty} (M)${\rm ,} and for some positive constant
$\delta${\rm ,} we have
\begin{equation}
\lim\limits_{r \rightarrow \infty}\frac{K(r)}{r^{1 + \delta}} = +
\infty,
\end{equation}
then there exists no complete metric $g$ on $M$ satisfying the
conditions {\rm (}a{\rm )}--{\rm (}c{\rm )} in Theorem~$3.6$.
\end{coro}

\begin{proof}
By (3.26), it is easy to show that there is $a > 0$ such that
$\int_{0}^{r}\frac{1}{V(t)}\int_{B(t)}K\d \mu \d t > 0$ for $r
\geq a$. To prove the corollary, it is sufficient to show that
condition (3.25) is satisfied. This can done by using the limit
comparison test for improper integral.

In fact, from the assumption of the theorem for sectional
curvature and the well-known Laplacian comparison theorem
\cite{6}, we have
\begin{equation}
\frac{n - 1}{r} \leq \Delta r \leq (n - 1)c\coth (cr).
\end{equation}

So $\Delta r$ is bounded as $r \rightarrow \infty$.

On the other hand, by (3.26) it is obvious that
\begin{equation*}
\lim\limits_{t\rightarrow \infty} \int_{B(t)} K\d \mu =
\lim\limits_{t\rightarrow \infty}\int_{0}^{t} K(s)V(s)\d s = +
\infty.
\end{equation*}

Thus by the L'H$\hat{\rm o}$pital's rule, (3.26) and (3.27), we
have
\begin{align*}
\lim\limits_{t\rightarrow \infty} \frac{\int_{B(t)}K\d \mu}{V(t)}
&= \lim\limits_{t\rightarrow \infty}\frac{K(t)V(t)}{V'(t)}\\[.3pc]
&= \lim\limits_{r\rightarrow \infty}\frac{K(r)}{\Delta r}\\[.3pc]
&\geq \lim\limits_{r\rightarrow \infty} \frac{K(r)}{(n -
1)c\coth(cr)}\\[.3pc]
&= + \infty.
\end{align*}

Then it is easy to see that
\begin{equation*}
\lim\limits_{r\rightarrow
\infty}\int_{0}^{r}\frac{1}{V(t)}\int_{B(t)} K \d \mu \d t = +
\infty.
\end{equation*}

Now by the L'H$\hat{\rm o}$pital's rule we have
\begin{align*}
\lim\limits_{r\rightarrow \infty} \frac{r^{2 +
\delta}}{\int_{0}^{r}\frac{1}{V(t)}\int_{B(t)}K\d \mu \d t} &=
\lim\limits_{r\rightarrow \infty} \frac{(2 + \delta)r^{1 +
\delta}V(r)}{\int_{0}^{r}K(t)V(t)\d t}\\[.3pc]
&= (2 + \delta) \lim\limits_{r\rightarrow \infty} \frac{(1 +
\delta)r^{\delta}V(r) + r^{1 + \delta}V'(r)}{K(r)V(r)}\\[.3pc]
&= (2 + \delta) \lim\limits_{r\rightarrow \infty} \frac{r^{1 +
\delta}}{K(r)} \left(\frac{1 + \delta}{r} + \Delta r\right)\\[.3pc]
&= 0.
\end{align*}

Therefore we obtain
\begin{equation*}
\lim\limits_{r\rightarrow \infty} \frac{r^{\frac{2 +
\delta}{2}}}{\left[\int_{0}^{r}\frac{1}{V(t)}\int_{B(t)}K{\rm d
}\mu{\rm d}t\right]^{1/2}} = 0.
\end{equation*}

Then (3.25) follows by the limit comparison test. This completes
the proof of the corollary. \hfill $\Box$
\end{proof}

\section*{Acknowledgement}

The author would like to thank Professor Xu Zongben for his
constant help. This work is supported by the Scientific Research Fund of
Wenzhou Normal College (No.~2003z12).

\end{document}